**L. Berezansky [a] and L. Idels [b]**

[a]Department of Mathematics and Computer Science, Ben-Gurion University of the Negev, Beer-Sheva 84105, Israel. Partially supported by the Israeli Ministry of Absorbtion.

[b] Author to whom all correspondence should be addressed. Mathematics Department, Malaspina University-College, Building 360/304, 900 Fifth Street, Nanaimo, British Columbia V9R 5S5, Canada. Partially supported by an NSERC Research Grant.

**Population Models With Delay in Dynamic Environment**

**Abstract**: We study the combined effects of periodically varying carrying capacity and survival rates on the fish population in the ocean (sea). We introduce the Getz type delay differential equation model with a control parameter which describes how fish are harvested. We will modify and extend harvesting model of an exploited fish population to include periodic and rotational harvesting rates.

We study the existence of global solutions for the initial value problem, extinction and persistence conditions, and the existence of periodic solutions.

**1. Introduction**

Traditional Population Dynamics is based on the concept that carrying capacity does not change over time, even though it is known that the values of carrying capacity related to the habitat areas might vary, e.g. the more biodiverse a system, the greater the carrying capacity and hence also the productivity [14,16,17].

Some scientists have recognized carrying capacity as problematic, and it is not a species concept, it is a system concept.

We will modify the Getz type difference equation model [8] and extend it to include varying carrying capacity, fecundity and natural mortality rates and periodic harvesting rate. We study the existence of positive global solutions for the initial value problem, persistence conditions, and the existence of periodic solutions.

**2. Preliminaries**

Consider the following standard differential equation which is widely used in Population Dynamics [5,7,15]

$$\dot{N} = [\beta(t,N) - \eta(t,N)]N - \lambda(t)N. \quad (1)$$

where $N = N(t)$ is the population biomass, $\beta(t,N)$ is the per-capita fecundity rate (the birth rate), $\eta(t,N)$ is the per-capita mortality rate, and $\lambda(t)$ is the harvesting rate per-capita.

A classical logistic harvesting model [7] has the following form

Here a function of captures is defined as $Y(t) = qN(t)E$, where $q \geq 0$ is the catchability coefficient, defined as the fraction of the population fished by an effort unit. Here $E \geq 0$ is the fishing effort.

It is a well-known fact [5,15] that the traditional logistic model in some cases produces artificially complex dynamics, therefore it would be reasonable to get away from the specific logistic form in studying population dynamics and use more general classes of growth models [1-4]. For example, in order to drop an unnatural symmetry of the logistic curve, recently we have considered [2], the modified delay logistic form of Pella and Tomlinson or Richards' growth delay equation

$$\dot{N}(t) = rN(t)\left[1 - \left(\frac{N(t-\tau)}{K}\right)^m\right].$$

In equation (1) let $\beta(t,N)$ be a Hill's type function

$$\beta(t,N) = \frac{r}{1 + (N/K)^\gamma},$$

where $r > 0$, $K > 0$. Parameter $\gamma > 0$ is referred to by W. Getz [8] as the "abruptness" parameter.

Generally, models with the delay in the reproduction term recognize that for real organisms it takes time to develop from newborns to reproductively active adults. If we take into account that delay and assume that $\eta(t,N) = \eta(t)$, then we have the following time-lag model based on equation (1)

$$\dot{N}(t) = \left[ \frac{r(t)}{1 + \left( \dfrac{N(g(t))}{K(t)} \right)^{\gamma}} - \eta(t) \right] N(t) - \lambda(t) N(t), \quad t \geq 0. \quad (2)$$

Here $r(t) > 0$ is a fecundity factor, $\eta(t) > 0$ is a mortality factor, $\lambda(t) > 0$ is a harvesting factor, $K(t) > 0$ is a carrying capacity, $g(t)$ is the time to develop from newborns to reproductively active adults, $0 \leq g(t) \leq t$. For $\lambda(t) \equiv 0$, $r(t) \equiv r$, $g(t) = t - \theta$, $\eta(t) \equiv \eta$, $K(t) \equiv K$, equation (2) is the Getz type differential equation [8] that describes dynamics of the fish population.

## 3. Proportional Harvesting in Dynamic Environment

If we denote $b(t) = \eta(t) - \lambda(t)$, then equation (2) has the following form

$$\dot{N}(t) = \frac{r(t) N(t)}{1 + \left[ \dfrac{N(g(t))}{K(t)} \right]^{\gamma}} - b(t) N(t), \quad t \geq 0, \quad (3)$$

with the initial function and the initial value

$$N(t) = \varphi(t), \, t < 0, \, N(0) = N_0, \quad (4)$$

under the following conditions:

(a1) $\gamma > 0$;

(a2) $r(t)$, $b(t)$, $K(t)$ are continuous on $[0. \infty)$ functions, $r(t) > 0$, $b(t) \geq b > 0$,

$K(t) \geq k > 0$;

(a3) $g(t)$ is a continuous function, $g(t) \leq t$; $\limsup\limits_{t \to \infty} g(t) = \infty$;

(a4) $\varphi : (-\infty, 0) \to R$ is a continuous bounded function, $\varphi(t) \geq 0, N_0 > 0$.

**Definition**. A function $N : R \to R$ with continuous derivative is called *a solution of problem* (3)-(4), if it satisfies equation (3) for all $t$ $[0,\infty)$ and equalities (4) for $t \leq 0$. If $t_0$ is the first point, where the solution $N(t)$ of (3)-(4) vanishes, i.e., $N(t_0) = 0$, then we consider the only positive solutions of the problem (3)-(4) on the interval $[0, t_0)$.

**Theorem 1.** Suppose (a1)-(a4) hold and

$$r(t) > b(t), \quad \sup_{t \geq 0} \int_{g(t)}^{t} [r(s) - b(s)]ds < \infty, \quad \sup_{t \geq 0} \int_{g(t)}^{t} b(s)ds < \infty.$$

Then problem (3)-(4) has on $[0, \infty)$ a unique positive solution $N(t)$ such that

$$\min\left\{N_0, \inf_{t \geq 0} K(t)\left(\frac{r(t)}{b(t)} - 1\right)^{\frac{1}{\gamma}} \exp\left\{-\sup_{t \geq 0} \int_{g(t)}^{t} b(s)ds\right\}\right\} \leq N(t) \tag{5a}$$

and

$$N(t) \leq \max\left\{N_0, \sup_{t \geq 0} K(t)\left(\frac{r(t)}{b(t)} - 1\right)^{\frac{1}{\gamma}} \exp\left\{\sup_{t \geq 0} \int_{g(t)}^{t} r(s) - b(s))ds\right\}\right\} \tag{5b}$$

**Proof**. The existence of the unique local solution is a consequence of well-known results for nonlinear delay differential equations (see, for example, [9-11]).

Clear

$$N(t) = N_0 \exp\left\{\int_0^t \left[\frac{r(s)}{1 + \left[\frac{N(g(s))}{K(s)}\right]^{\gamma}} - b(s)\right] ds\right\}$$

hence the local solution is positive.

If $[0,\alpha)$ is the maximal interval of the existence of this solution, where $\alpha < \infty$, then

$\lim_{t \to \alpha^-} N(t) = +\infty.$

Eq.(3) implies that $\dot{N}(t) \leq [r(t) - b(t)]N(t),$ therefore

$$N(t) \leq \exp\left\{\int_0^t [r(s) - b(s)]ds\right\}.$$

Then the local solution $N(t)$ is bounded on the maximal interval $[0,\alpha)$ of the existence of this solution, and we have a contradiction. Hence the maximal interval of the existence of $N(t)$ is $[0,\infty)$ and the global solution is positive. Therefore we have to prove only inequalities (5a-5b).

Suppose that $\dot{N}(t) > 0$ for any $t > 0$. Then (3) implies that

$$\frac{r(t)}{1 + \left[\dfrac{N(g(t))}{K(t)}\right]^\gamma} - b(t) > 0$$

or $\quad N(g(t)) < K(t)\left(\dfrac{r(t)}{b(t)} - 1\right)^{1/\gamma}$

and we have $\dot{N}(t) \leq N(t)(r(t) - b(t))$.

Finally

$$N(t) \leq N(g(t))\exp\left\{\int_{g(t)}^t [r(s) - b(s)]ds\right\} \leq K(t)\left(\frac{r(t)}{b(t)} - 1\right)^{1/\gamma} \exp\left\{\int_{g(t)}^t [r(s) - b(s)]ds\right\}$$

.

The last inequality proves (5b) in the case where $\dot{N}(t) > 0$.

If $\dot{N}(t) < 0$ for any $t > 0$ inequality (5b) is evident.

Suppose now that $\{t_n\}$ is the sequence of the points where the function $N(t)$ has the local maximum. Then $\dot{N}(t_n) = 0$, hence

$$N(g(t_n)) = K(t_n)\left(\frac{r(t_n)}{b(t_n)} - 1\right)^{1/\gamma}.$$

In the interval $[g(t_n), t_n]$ we have $\dot{N}(t) \leq N(t)(r(t) - b(t))$, then

$$N(t) \leq N(g(t_n))\exp\int_{g(t_n)}^{t}[r(s) - b(s)]ds.$$

Hence $N(t_n) \leq N(g(t_n))\exp\left\{\int_{g(t_n)}^{t_n}[r(s) - b(s)]ds\right\}.$

Inequality (5b) is satisfied since $\sup_{t \geq 0} N(t) = \sup_n N(t_n)$.

Suppose now that $\dot{N}(t) < 0$ for any $t > 0$. Then

$$\frac{r(t)}{1 + \left[\frac{N(g(t))}{K(t)}\right]^{\gamma}} - b(t) > 0$$

and

$$N(g(t)) > K(t)\left(\frac{r(t)}{b(t)} - 1\right)^{1/\gamma}.$$

Equation (3) implies $\dot{N}(t) \geq -b(t)N(t)$. Therefore

$$N(t) \geq N(g(t))\exp\left\{-\int_{g(t)}^{t}b(s)ds\right\} \geq K(t)\left(\frac{r(t)}{b(t)} - 1\right)^{1/\gamma}\exp\left\{-\int_{g(t)}^{t}b(s)ds\right\}.$$

Hence inequality (5a) has been proven for $\dot{N}(t) < 0$.

If $\dot{N}(t) > 0$ for any $t > 0$, then inequality (5a) is evident.

Suppose now that $\{\tau_n\}$ is the sequence of the points where the function $N(t)$ has the local minimum. Then $\dot{N}(\tau_n) = 0$ and

$$N(g(\tau_n)) = K(\tau_n)\left(\frac{r(\tau_n)}{b(\tau_n)} - 1\right)^{1/\gamma}.$$

In the interval $[g(\tau_n), \tau_n]$ we have $\dot{N}(t) \geq -b(t)N(t)$, therefore

$$N(t) \geq N(g(\tau_n))\exp\left\{-\int_{g(t_n)}^{t} b(s)ds\right\}.$$

and

$$N(\tau_n) \geq N(g(\tau_n))\exp\left\{-\int_{g(t_n)}^{\tau_n} b(s)ds\right\}.$$

Since $\inf_{t \geq 0} N(t) = \inf_n N(\tau_n)$ inequality (5a) is true in this case, therefore we have proved Theorem 1.

Remark. The inequalities (5a -5b) imply that the solution $N(t)$ of problem (3)-(4) is persistent.

## 4. Periodic Proportional Harvesting in a Dynamic Environment

Periodic harvesting is used frequently as a tool by fishery managers to protect a stock during a spawning season for a defined period of time [5-7].

Consider again equation (2).

$$\dot{N}(t) = \frac{r(t)N(t)}{1 + \left[\dfrac{N(g(t))}{K(t)}\right]^{\gamma}} - \eta(t)N(t) - \lambda(t)N(t). \tag{6}$$

Now we assume that *r, K, η* and *λ* are all positive continuous T-periodic functions.

Function *g(t)* is defined as $g(t) = t - \theta(t)$ where $\theta(t)$ is a T-periodic function.

If fish reproduce primarily in the spring of the year, the functions *r(t)* and *η(t)* might be modeled by some periodic functions, such as

$r(t) = r_0 + A\cos(\varpi\pi(t - 0.25))$, where $0 < A \leq r_0$,

$\eta(t) = \eta_0 + \cos(\varpi\pi(t - 0.25))$, where $0 < \eta \leq \eta_0$.

If the food supply peaks each year in the fall, then $K(t) = K_0 + C\cos(\varpi\pi(t - 0.75))$, where $0 < C \leq K_0$.

To introduce a periodic harvesting rate one might use for example

$$\lambda(t) = \begin{cases} 0.5\sin\dfrac{\pi(t - n - t_{start})}{H} & \text{if } n + t_{start} < t < n + t_{start} + H, \quad n = 0,1,2,..., \\ 0 & \text{otherwise,} \end{cases}$$

where *H* is the harvesting time, $t_{start}$ is the harvest starting time within one year, *n* is the number of years, e.g., if $H = 0.25$ and $t_{start} = 0.25$ (harvest in the summer season only).

There has been growing interest in rotational use of marine resource areas [13]. Rotational fishing is usually thought of as a sequence of periodic closure and openings of different fishing areas. For example, three-year pulse rotation, when the area is closed for two years and fished for one year.

To model a rotational harvest, say, three-year rotation, one could use

$$\lambda(t) = \begin{cases} 0.5k\sin\dfrac{\pi(t - n - t_{start})}{H} & \text{if } n + t < t < n + t_{start} + H, \quad k = [n/3], n = 0,1,2,... \\ 0 & \text{other} \end{cases}$$

Denote $b(t) = \eta(t) - \lambda(t)$ and assume that $g(t) = t - \theta(t)$.

Equation (6) takes the following form

$$\dot{N}(t) = \left[ \frac{r(t)}{1 + \left[ \frac{N(t-\theta(t))}{K(t)} \right]^{\gamma}} - b(t) \right] N(t) \qquad (7)$$

**Theorem 2**. Let *b(t), K(t), θ(t)* and *r(t)* be T-periodic functions. Suppose *r(t) > b(t)*.

If at least one of the following conditions hold:

(b1) $\quad m = \inf\limits_{t \geq 0} \left( \frac{r(t)}{b(t)} - 1 \right) K^{\gamma}(t) > 1$

(b2) $\quad M = \sup\limits_{t \geq 0} \left( \frac{r(t)}{b(t)} - 1 \right) K^{\gamma}(t) < 1$

then equation (7) has at least one periodic positive solution $N_0(t)$.

To prove that result we will use the following Lemma [13].

**Lemma 1**. Consider the delay differential equation

$$\frac{dN}{dt} = \pm N(t) G(t, N(t - \theta(t, N(t)))) \qquad (8)$$

where $G, \theta \in C(R^2, R)$ and $G, \theta$ are T-periodic functions with respect to the first argument.

Suppose that there exists constants $B, \alpha, \beta > 0$ such that:

(I)      when $|x| < B$, the inequality $|G(t, e^x)| \leq \beta$ holds uniformly for $t \in R$,

and

        when $|x| > B$, the inequality $xG(t, e^x) > 0$ holds uniformly for $t \in R$.

(II) one of the following conditions holds:

     (i) when $x < -B$, $G(t, e^x) > -\alpha$ holds uniformly for $t \in R$.

     (ii) when $x > B$, $G(t, e^x) < \alpha$ holds uniformly for $t \in R$.

Then equation (8) has at least one positive T-periodic solution.

**Proof of Theorem 2.** Denote

$$G(t,u) = b(t) - \frac{r(t)}{1 + \left[\dfrac{u}{K(t)}\right]^{\gamma}}.$$

then

$$G(t,e^x) = b(t) - \frac{r(t)}{1 + \left[\dfrac{e^x}{K(t)}\right]^{\gamma}}$$

We have

$$\left|G(t,e^x)\right| \leq b(t) + r(t) \leq \sup_{t \geq 0}[b(t) + r(t)] = \beta.$$

Hence for every $B > 0$ the first part of statement (I) of Lemma 1 holds for $|x| < B$.

Now we have to prove that for $|x| > B$ the inequality $xG(t,e^x) > 0$ is true, where the number $B$ will be defined explicitly.

Suppose (b1) holds. The inequality $G(t,e^x) > 0$ is equivalent to

$$x > \frac{1}{\gamma} \ln\left\{\left(\frac{r(t)}{b(t)} - 1\right)K^{\gamma}(t)\right\}$$

Denote

$$B = \frac{1}{\gamma} \sup_{t \geq 0} \ln\left\{\left(\frac{r(t)}{b(t)} - 1\right)K^{\gamma}(t)\right\}$$

Condition (b1) implies that $B > 0$.

Hence for $x > B$ we have $xG(t,e^x) > 0$.

Let $x < -B$. The inequality $G(t,e^x) < 0$ is equivalent to

$$x < \frac{1}{\gamma} \ln\left\{\left(\frac{r(t)}{b(t)} - 1\right) K^{\gamma}(t)\right\}.$$

This inequality is satisfied by (b1):

$$\ln\left\{\left(\frac{r(t)}{b(t)} - 1\right) K^{\gamma}(t)\right\} > 0.$$

Therefore it follows that $xG(t,e^x) > 0$.

To check the first condition of part (II) of Lemma 1 we have

$$G(t,e^x) = b(t) - \frac{r(t)}{1 + \left[\dfrac{e^x}{K(t)}\right]^{\gamma}} \geq b(t) - r(t) > -\alpha.$$

Hence for every $x$ we have $G(t,e^x) > -\alpha$. Then first part of condition (II) of Lemma 1 holds, therefore equation (7) has at least one positive T-periodic solution.

If the second condition ($b2$) is satisfied then the proof of the Theorem 2 is similar.

## 5. References.